\documentclass[12pt]{article}

\usepackage{comment}
\usepackage{amssymb}
\usepackage{amsfonts}
\usepackage{mathrsfs}
\usepackage{enumerate}
\usepackage{amssymb}
\usepackage{xcolor}
\usepackage{color}
\usepackage{colortbl}
\usepackage{bm}
\usepackage{natbib}
\usepackage{graphics}
\usepackage{booktabs}
\usepackage{amsfonts}
\usepackage[T1]{fontenc}
\usepackage{multirow}
\usepackage[font=small,labelfont=bf,tableposition=top]{caption}
\usepackage[subrefformat=parens]{subcaption}
\usepackage{threeparttable}
\usepackage{booktabs}
\usepackage{array}
\usepackage[dvips]{graphicx}
\usepackage{amsmath}
\DeclareMathOperator*{\argmax}{argmax}
\usepackage{enumerate}
\usepackage{natbib,epsfig,graphicx}
\usepackage{url}
\usepackage{epstopdf}
\usepackage{algorithm}
\usepackage{algorithmic}
\usepackage{mathrsfs}
\usepackage{amsmath}
\usepackage{mathrsfs}
\usepackage{dsfont}
\usepackage{mathrsfs}
\usepackage{JASA_manu}

\makeatletter
\def\singlespace{\def\baselinestretch{1}\@normalsize}

\newtheorem{definition}{Definition}[subsection]

\newtheorem{proposition.A}{Proposition}[section]

\newtheorem{theorem}{Theorem}[subsection]

\newtheorem{remark.A}{Remark}[section]
\newtheorem{assumption}{Assumption}[section]
\begin{document}

\renewcommand{\baselinestretch}{1.5}

\title{Hierarchical Change-Point Detection for Multivariate Time Series via a Ball Detection Function}
\author{{\sc Xueqin Wang, Qiang Zhang, Wenliang Pan, Xin Chen}\\{\sc and Heping Zhang}\\
{Sun Yat-Sen University}\\{Southern University of Science and Technology}\\{Yale University}}

\date{}
\maketitle

\begin{singlespace}
\begin{footnotetext}[1]
{Xueqin Wang is Professor, Department of Statistical Science, School of
	Mathematics, Southern China Center for Statistical
	Science, Sun Yat-Sen University, Guangzhou, 510275, China; Zhongshan School of Medicine, Sun Yat-Sen
	University, Guangzhou, 510080, China; and Xinhua College, Sun Yat-Sen University, Guangzhou, 510520, China (e-mail:
	wangxq88@mail.sysu.edu.cn).
Qiang Zhang is Ph.D. candidate, Department
of Statistical Science, School of Mathematics, Sun Yat-Sen University, Guangzhou, 510275, China (e-mail:zhangqg2@mail2.sysu.edu.cn). Wenliang Pan is research associate, Department
of Statistical Science, School of Mathematics, Sun Yat-Sen University, Guangzhou, 510275, China (e-mail:panwliang@mail.sysu.edu.cn).
Xin Chen is associate professor, Department of Mathematics,  Southern University of Science and Technology (e-mail: chenx8@sustc.edu.cn).
Heping Zhang is Susan Dwight Bliss Professor of Biostatistics, Yale
University School of Public Health, New Haven, CT 06520 (e-mail: heping.zhang@yale.edu).
Zhang's research was partially supported by grant R01 MH116527 from the National Institute of Mental Health and DMS-1722544 from the National Science Foundation.}

\end{footnotetext}
\end{singlespace}

\newpage
\begin{center}
\bf {Hierarchical Change-Point Detection for Multivariate Time Series via  a Ball Detection Function}
\end{center}

\begin{abstract}

\baselineskip=20pt
 Sequences of random objects arise from many real applications, including high throughput omic data and  functional imaging data. Those sequences are usually dependent, non-linear, or even Non-Euclidean, and an important problem is  change-point detection in such dependent sequences in Banach spaces or metric spaces. The problem usually  requires the accurate inference for not only whether changes might have occurred but also the locations of the changes when they did occur. To this end, we first introduce a Ball detection function and show that it reaches its maximum at the change-point if a sequence has only one change point. Furthermore, we propose a  consistent estimator of  Ball detection function based on which we develop a hierarchical algorithm to detect all possible change points. We prove that the estimated change-point locations are consistent.  Our procedure can  estimate the number of change-points and detect their locations without assuming any particular types of change-points as a change can occur in a sequence in different ways.  Extensive simulation studies and analyses of two interesting real datasets wind direction and Bitcoin price  demonstrate that our method has considerable advantages over existing competitors, especially when data are non-Euclidean or when there are distributional changes in the variance.
\end{abstract}
\bigskip

\noindent\textsc{\bf Keywords}: {Ball Divergence; Change point detection;  Non-Euclidean data; Ergodic stationary sequence; Absolutely regular sequence}.
 \baselineskip=22pt

\newpage

\section{Introduction}
 Stationarity is crucial in analyzing random sequences because statistical inference usually requires a probabilistic mechanism constant in, at least, a segment of observations.  Therefore, it is important to detect whether changes occur in a sequence of observations prior to statistical inference.  Such change-point problems arise from many applications: abrupt events in video surveillance \citep{mayer2015change};   deterioration of product quality in quality control \citep{lai1995sequential}; credit card fraud in finance \citep{bolton2002statistical}; alteration of genetic regions in cancer research \citep{erdman2008fast} and so on.

There is a large and rapidly growing literature on change-point detection \citep{aminikhanghahi2016survey,niu2016multiple,sharma2016trend,fryzlewicz2014wild}. Many methods rely on the assumed parametric models to detect special change types  such as location, scale or presumed distribution family. \citet{page1954continuous} introduced a method by examining the ratio of log-likelihood functions. \citet{lavielle2006detection} detected change-points by maximizing a log-likelihood function. \citet{yau2016inference} proposed a likelihood ratio scan method for piecewise stationary autoregressive time series. Some Bayesian change-point detection methods assume that the observations are normally distributed, and calculate the probability of change-point at each point \citep{barry1993bayesian,Zhang1995changepoint,wang2015bayesian,maheu2018efficient} to name a few. Since parametric methods potentially suffer from model misspecification, other methods are developed to detect general distributional changes with more relaxed assumptions.
\citet{kawahara2012sequential} provided an algorithm which relied heavily on estimating the ratio of probability densities. \citet{lung2015homogeneity} identified change-points via the well-known Wilcoxon rank statistic. \citet{matteson2014nonparametric} proposed a nonparametric method using the concept of energy distance for independent observations. There are also some binary segmentation methods statistics  \citep{fryzlewicz2014wild,cho2015multiple,eichinger2018mosum}. Two advantages of the binary segmentation procedures are their simplicity and computational efficiency, but their false discovery rates may be hard to control because they are `greedy' procedures. \citet{zou2014nonparametric} introduced a nonparametric empirical likelihood approach to detecting multiple change-points in independent sequences, and estimated the locations of the change-points by using the dynamic programming algorithm and the intrinsic order structure of the likelihood function.

Automatically detecting the number of change-points is also important.
Some methods are developed to detect only a single change-point \citep{ryabko2008hypotheses}, while some methods require a known number of change-points but unknown locations \citep{hawkins2001fitting,lung2015homogeneity}. In real data analysis, however, we usually do not know the number of change-points.

With increasing richness of data types,  non-Euclidean data, such as shape data, functional data, and spatial data, commonly arise from applications. For example, one of the problems of interest to us is the changes in the monsoon direction as defined by circle, a simple Riemannian manifold. Methods developed in Hilbert spaces are not effective for this type of problems as our analysis of the data from Yunnan-Guizhou Plateau ($105^\circ$E, $27^\circ$N) collected from 2015/06/01 to 2015/10/30 illustrates below.  To the best of our knowledge, few methods exist to detect change-points in a non-Euclidean sequence.  \citet{chen2015graph} and \citet{chu2019asymptotic} proposed a series of graph-based nonparametric approaches that could be applied to non-Euclidean data with arbitrary dimension. However, their proposed  methods apply to \emph{iid} observations only and are restricted to one or two change-points. Therefore, it remains to be an open and challenging problem to develop methods to detect arbitrarily distributional changes for non-Euclidean sequences, including the change-point locations and the number of the change-points.

To address this challenge,  we introduce a novel concept of Ball detection function via Ball divergence \citep{pan2018ball}. Ball divergence  is a recently developed measure of divergence between two probabilities in separable Banach spaces. The Ball divergence is zero if and only if the two probability measures are identical. Since its sample statistic is constructed by metric ranks, the test procedure for an identical distribution is robust to heavy-tailed data or outliers, consistent against alternative hypothesis, and applicable to imbalanced data. Therefore, the empirical ball divergence is an ideal statistic to test whether or not a change has occurred. Unfortunately, it does not inform us where the change occurs, because in theory the probability measures before and after any time point are always different if there exists a change point in the sequence. Therefore it is imperative for us to observe how the probability measures before and after any time vary with time and then develop a proper criterion to detect the change-point location. We introduce a Ball detection function as an effective choice which reaches its maximum at the change point if a sequence has only one change point. We further develop a hierarchical algorithm to detect multiple change-points using the statistic based on the Ball detection function. The advantages of our procedure are threefold: our procedure can  estimate the number of change-points and detect their locations;  our procedure can detect any types of change-points; and both uniquely and importantly, our procedure can handle complex stochastic sequences, for example, non-Euclidean sequences.

The rest of this article is organized as follows. In Section 2, we review the notion of Ball divergence,  and then introduce a novel change-point detection function, i.e., a Ball detection function based on Ball divergence with a scale parameter for weakly dependent sequences. We further establish its asymptotic properties. We show how to use the Ball detection function to detect change-points and  establish the consistent properties of our method in Section 3.  In Section 4, we compare the performance of our method with some existing methods in various simulation settings. In section 5, two real data analyses demonstrate the utility  of our proposed method. We make some concluding remarks in Section 6. All technical details are deferred to Appendix.

\section{Change-point Detection in Dependent Sequences}

\subsection{Review of Ball Divergence}

Ball divergence (BD, \citet{pan2018ball})  is a measure of the difference between two probabilities in a separable Banach space $(\mathrm{A},||\cdot||)$, with the norm $||\cdot||$. $\forall~u,v \in \mathrm{A}$, the distance  between $u$ and $v$ deduced from the norm is $\rho(u,v)=|| u-v ||$. Denote by $\bar{B}(u,r)=\{x|\rho(x,u)\leq r\}$  a closed ball. Let $\mathcal{B}$ be the smallest $\sigma$-algebra in $\mathrm{A}$ that contains all closed (or open) subsets of $\mathrm{A}$. Let $\mu$ and $\nu$ be two probabilities on $\mathcal{B}$. Ball divergence \citep{pan2018ball} is defined as follows.

\begin{definition}\label{df_bdo}
	The Ball divergence of  two Borel probabilities $\mu$ and $\nu$ in $\mathrm{A}$ is defined as an integral of the square of the measure difference between $\mu$ and $\nu$ over  arbitrary closed balls,
	\begin{equation*}%\label{BD}
	D(\mu,\nu)=\iint_{\mathrm{A}\times \mathrm{A}} [\mu-\nu]^2 (\bar{B}(u,\rho(u,v)))(\mu(du) \mu(dv)+ \nu(du)\nu(dv)).
	\end{equation*}
\end{definition}

Let $S_\mu$ and $S_\nu$ be the support sets of $\mu$ and $\nu$ respectively. The BD has the following important property \citep{pan2018ball}:

\begin{theorem}\label{quantity} Given two Borel probabilities $\mu$ and $\nu$ in a  finite dimensional Banach space $\mathrm{A}$, then $D(\mu,\nu)\geq 0$  where the equality holds if and only if $\mu=\nu$.  It can be extended to separable Banach spaces if $S_\mu=\mathrm{A}$ or $S_\nu=\mathrm{A}$.
\end{theorem}

\subsection{Ball Divergence with a Scale Parameter}

The Ball divergence introduced above cannot detect the locations of change-points accurately enough while comparing the distributions of the sequences before and after the change-points. We need to introduce a  Ball divergence associated with a scale parameter $\alpha$ as follows.

\begin{definition}\label{df_bd} A Ball divergence of two Borel measures $\mu$ and $\nu$ in $\mathrm{A}$ is defined as
	\begin{align}\label{formula_alphaBD}
	D_{\alpha}(\mu,\nu)
	=\iint_{\mathrm{A}\times \mathrm{A}} [\mu-\nu]^{2}(\bar{B}(u,\rho(u,v))\omega_{\alpha} (du)\omega_{\alpha} (dv),
	\end{align}
	where $\omega_{\alpha}=\alpha \mu+(1-\alpha) \nu$ is the mixture distribution measure with the scale parameter $\alpha\in[0,1]$.
\end{definition}

$D_{\alpha}(\mu,\nu)$ also has the equivalence property below, which is critical to the comparison of the distributions of any two sequences.
\begin{theorem}\label{BDifandonlyif}
	Given two Borel probabilities $\mu,\nu$ in a finite dimensional Banach space $\mathrm{A}$, then $D_{\alpha}(\mu,\nu)\geq 0$ where the equality holds if and only if $\mu=\nu$. It also holds on separable Banach spaces if $S_\mu=\mathrm{A}$ or $S_\nu=\mathrm{A}$.
\end{theorem}

Theorem \ref{BDifandonlyif} assures that for any $\alpha\in [0,1],$ $D_{\alpha}(\mu,\nu)$ possesses the most important property as $D(\mu,\nu)$ in terms of
testing the distributional difference between two sequences.  Importantly, with the introduction of $\alpha$, we can consistently estimate the locations of the change-points. Here, we highlight the relationship and difference between $D_{\alpha}(\mu,\nu)$ and $D(\mu,\nu).$

When $\alpha=1$, $D_{1}(\mu,\nu)$ is the measure difference over the balls whose centers and the endpoints of the radius following measure $\mu.$ When $\alpha=0$, $D_{0}(\mu,\nu)$ is the measure difference over the balls whose centers and the endpoints of the radius following the measure $\nu$. Moreover,
\begin{equation*}
D(\mu,\nu)=D_{0}(\mu,\nu)+D_{1}(\mu,\nu).
\end{equation*}
For $\alpha\in(0,1)$, $D_{\alpha}(\mu,\nu)$ is the mean of the measure differences from two samples over the balls whose centers and endpoints of the radius following four possible pairs of measures:$(\mu,\mu)$, $(\mu,\nu)$,$(\nu,\mu)$, and $(\nu,\nu)$ where the ratio of two measures is $\alpha:1-\alpha$.

Ball divergence with a  scale parameter can be defined in the general metric space, following the Generalized Banach-Mazur theorem \citep{kleiber1969a} as stated in the Supplementary material.

\subsection{Ball Detection Function}
Now, we introduce a Ball detection function which is maximized at the change point if there exists one, and hence can be used to
determine the location of the change point. For clarity, let us consider a conceptual sequence with a change point $\alpha \in (0,1)$, and the probability measures before and after $\alpha$ are $\mu$ and $\nu$, respectively.  Denote the indicator function by $I(\cdot)$. For  a "time" $\beta \in (0,1)$, define 
\begin{equation*}
  h_{\alpha}(\beta)=\frac{\alpha}{\beta}I(\beta\geq\alpha)+\frac{1-\alpha}{1-\beta}I(\beta<\alpha).
\end{equation*}
Without loss of generality, suppose that $\beta>\alpha$, the probability measures before and after $\beta$ are $\frac{\alpha}{\beta}\mu+(1-\frac{\alpha}{\beta})\nu$ and $\nu$. By the definition of Ball divergence (\ref{formula_alphaBD}), we have
\begin{equation*}
D_{\beta}(\frac{\alpha}{\beta}\mu+(1-\frac{\alpha}{\beta})\nu,\nu)=(\frac{\alpha}{\beta})^2D_{\alpha}(\mu,\nu).
\end{equation*}
Therefore, in general,
\begin{equation}\label{relation alpha-beta}
D_{\beta}(h_{\alpha}(\beta)\mu+(1-h_{\alpha}(\beta)
)\nu,\nu)=h^2_{\alpha}(\beta)D_{\alpha}(\mu,\nu).
\end{equation}
The maximum of $h_{\alpha}(\beta)$ is attained when $\beta=\alpha$ if there exists a change-point $\alpha.$ In this case, we can find the change point by maximizing the ball divergence in equation (\ref{relation alpha-beta}). But we still need to test whether a change point has occurred or not. Next, we introduce a Ball detection function to simultaneously test the existence of a change-point and determine its location:
\begin{eqnarray*}\label{relation Ball-detection}
V(\beta;\mu,\nu) &=&\beta(1-\beta)D_{\beta}(h_{\alpha}(\beta)\mu+(1-h_{\alpha}(\beta)
)\nu,\nu)\\&=&\beta(1-\beta)h^2_{\alpha}(\beta)D_{\alpha}(\mu,\nu).
\end{eqnarray*}
Note that the maximum of $\beta(1-\beta)h_{\alpha}(\beta)$ is also attained when $\beta=\alpha$, allowing us to find the change point by maximizing $V(\beta;\mu,\nu)$. In next subsection, we shall discuss how this function is used to construct a test for a change-point test statistic.

\subsection{Ball Detection Function in Sample}
Suppose that a sequence of observations $\{Z_i\}_{1 \leq i \leq T}$ is comprised of  two multivariate stationary sequences $\{Z_i\}_{1 \leq i \leq M}$ with the probability measure $\mu_1$ and $\{Z_i\}_{M+1\leq i \leq T}$ with $\mu_2$, where  both  $\mu_1$ and $\mu_2$ are unknown. We estimate $D_{\alpha}(\mu_1,\mu_2)$ with $\alpha=M/T$ based on $\{Z_i\}_{1 \leq i \leq T}$.  Let $c(x,y;z)=I(z\in \bar{B}(x,\rho (x,y)))$, which identifies whether the point $z$  falls into the closed ball $\bar{B}(x,\rho (x,y))$ with $x$ as the center and $\rho (x,y)$ as the radius, and $e(x,y,z_1,z_2)=c(x,y;z_1)c(x,y;z_2)$, which determines whether two points $z_1$ and $z_2$ fall into the ball $\bar{B}(x,\rho (x,y))$ together.   Let $N=T-M$, $C_{ij}^{1}=\frac{1}{M}\sum_{u=1}^Mc(Z_i,Z_j;Z_{u}),
C_{ij}^{2}=\frac{1}{N}\sum_{v=M+1}^Tc(Z_i,Z_j;Z_{v}).$
A consistent estimator of the Ball divergence of $\mu_1$ and $\mu_2$  with the scale parameter $\alpha$ is
\begin{equation*}
D_{M,N}=\frac{1}{T^2}\sum\limits_{i,j=1}^{T}(C_{ij}^{1}-C_{ij}^{2})^2,
\end{equation*}
as summarized in Theorem \ref{slln}.

We also prove that $\frac{MN}{T}D_{M,N}$ has a limiting distribution under the null hypothesis in Theorem \ref{h0}. For this reason, we choose
	\begin{equation*}
	V(M,T)=\frac{MN}{T}D_{M,N}
	\end{equation*}
as the statistic to detect change-points.

To investigate the asymptotic properties of $V(M, T)$, we introduce two concepts of the random sequence: absolutely regular and ergodic stationary sequence.

Given the probability space $(\Omega,\mathscr{F},P)$ and two sub-$\sigma$-fields $\mathscr{A}$ and $\mathscr{B}$ of $\mathscr{F}$, let
$$
\beta(\mathscr{A},\mathscr{B})=\sup\sum_{i=1}^m\sum_{j=1}^n|P(A_{i}\bigcap B_{j})-P(A_{i})P(B_{j}))|,
$$
where the supreme is taken over all partitions of $\Omega$ into sets $A_1,\ldots,A_m\in\mathscr{A}$, all partitions of $\Omega$ into sets $B_1,\ldots,B_n\in\mathscr{B}$ and all $m,n\geq1$. A stochastic sequence $\{Z_i\}_{i\in\mathds{Z}}$ is called absolutely regular ( \citep{dehling2012asymptotic}, also called weakly Bernoulli \citep{aaronson1996strong}), if
$$
\beta(l)=\sup\limits_{n}\beta(\mathscr{F}_{0}^n,\mathscr{F}_{n+l}^\infty)\rightarrow 0,
$$
as $l\rightarrow\infty$. Here the $\mathscr{F}_i^j$ denotes the $\sigma$-field generated by the random variables $Z_i,\ldots,Z_j$. In this paper, we suppose that $\beta(l)=O(l^{-1-r})$ for any $r>0$. The concept of absolutely regular sequence is wide enough to cover all relevant examples from statistics except for long memory sequences.

Recall that an ergodic, stationary sequence (ESS) \citep{aaronson1996strong} is a random sequence $\{Z_i\}_{1\leq i\leq T}$ of form $Z_{i}=f(G^i)$ where $G^i$ is an ergodic, probability-preserving transformation in the probability space $(\Omega,\mathscr{F},P)$, and $f$ is a measurable function.
In essence, an ESS implies that the random sequence will not change its statistical properties with time (stationarity) and that its statistical properties  can be deduced from a single, sufficiently long sample of the sequence (ergodicity).

We have the following theorem for an absolutely regular sequence comprised of two ergodic stationary sequences:

\begin{theorem}\label{slln} Suppose that $\{Z_i\}_{1 \leq i \leq T}$ is an absolutely regular sequence, $\{Z_i\}_{1 \leq i \leq M}$ and $\{Z_i\}_{M+1\leq i \leq T}$ are both ergodic stationary with marginal probability measure $\mu_1,\mu_2$ respectively. When $M,T\rightarrow \infty$, $M/T\rightarrow\alpha_1$ for some $\alpha_1\in[0,1]$,  then
	$$\frac{V(M,T)}{T}\xrightarrow[M,T\rightarrow\infty]{a.s.}V(\alpha_1;\mu_1,\mu_2).$$
\end{theorem}

Theorem \ref{slln} means that $\frac{V(M,T)}{T}$ converges to Ball detection function $V(\alpha_1;\mu_1,\mu_2)$ almost surely. We further investigate the asymptotic distribution of $V(M,T)$. Under the null hypothesis, the Ball detection function in sample is the sum of four degenerate V-statistics. As in \cite{pan2018ball}, we denote $Q(x,y;x',y')$ as the second component in the H-decomposition of $V(M,T)$. Then we have the spectral decomposition:
$$Q(x,y;x',y')=\sum_{k=1}^{\infty}\lambda_kf_k(x,y)f_k(x',y'),$$
where $\lambda_k$ and $f_k$ are the eigenvalues and eigenfunctions of $Q(x,y;x',y')$. Let $\{Z_i'\}_{1\leq i\leq T}$ be an independent copy of $\{Z_i\}_{1\leq i\leq T}$. For $k\in\{1,2,\ldots\},$    $N_{1k},N_{2k}$ are assumed to be \emph{iid} $N(0,1)$, and let
\begin{align*}
a_k^2(\alpha_1)=(1-\alpha_1)E_{Z_1}[E_{Z_1'}f_k(Z_1,Z_1')]^2,\quad
b_k^2(\alpha_1)=\alpha_1E_{Z_1'}[E_{Z_1}f_k(Z_1,Z_1')]^2,\\
c_k^2(\alpha_1)=a_k^2(\alpha_1)+2(1-\alpha_1)(\sum_{j=1}^{\infty}E_{Z_1,Z_{1+j}}[E_{Z_1'}f_k(Z_1,Z_1')E_{Z_1'}f_k(Z_{1+j},Z_1')]),\\
d_k^2(\alpha_1)=b_k^2(\alpha_1)+2\alpha_1(\sum_{j=1}^{\infty}E_{Z_1',Z_{1+j}'}[E_{Z_1}f_k(Z_1,Z_1')E_{Z_1}f_k(Z_1,Z_{1+j}')]),
\end{align*}
$$\theta=E[E(c(Z_1,Z_2,Z_i)(1-c(Z_1,Z_2,Z_{j}))|Z_1,Z_2)].
$$
\begin{theorem}\label{h0}
	Under null hypothesis $H_{0}: \mu_1=\mu_2$, $\{Z_i\}_{1\leq i\leq T}$ is a stationary absolutely regular sequence with coefficients satisfying $\beta(l)=O(l^{-1-r})$ for $r>0$, if  $M,T\rightarrow \infty$, $M/T\rightarrow\alpha_1$ for some $\alpha_1\in[0,1]$, we have 
$$V(M,T)\xrightarrow[M,T\rightarrow\infty]{d}\sum_{k=1}^{\infty}\lambda_{k}[(c_{k}(\alpha_1)N_{1k}+d_{k}(\alpha_1)N_{2k})^{2}-(a_{k}^{2}(\alpha_1)+b_{k}^{2}(\alpha_1))]+\theta. $$
\end{theorem}

Under the alternative hypothesis, the Ball detection function in sample  is asymptotically normal because it is a sum of non-degenerate V-statistics. Let $g^{(1,0)}(Z_\mu)$ and $g^{(0,1)}(Z_{\nu})$ be the first component in H-decomposition of $V(M,T)$ and
$$\delta_{1,0}^2=Var(g^{(1,0)}(Z_u))+2\sum_{i=1}^\infty Cov(g^{(1,0)}(Z_u),g^{(1,0)}(Z_{u+i})),$$
$$\delta_{0,1}^2=Var(g^{(0,1)}(Z_v))+2\sum_{i=1}^\infty Cov(g^{(0,1)}(Z_v),g^{(0,1)}(Z_{v+i})).$$
% $\delta_{1,0}^2=Var(g^{(1,0)}(Z_\mu))$ and $\delta_{0,1}^2=Var(g^{(0,1)}(Z_{\nu}))$.
 We can obtain the asymptotic distribution under the alternative hypothesis.
\begin{theorem}\label{h1}
	$\{Z_i\}_{1\leq i\leq T}$ is a absolutely regular sequence with coefficients satisfying $\beta(l)=O(l^{-1-r})$ for $r>0$. Under $H_{1}: \mu_1\neq  \mu_2$, if $M,T\rightarrow \infty$, and  $M/T\rightarrow\alpha_1$ for some $\alpha_1\in[0,1]$, then we have
	$$\sqrt{\frac{T}{MN}}(V(M,T)-TV(\alpha_1;\mu_1,\mu_2))\xrightarrow[M,T\rightarrow\infty]{d} N(0,(1-\alpha_1)\delta_{1,0}^2+\alpha_1\delta_{0,1}^2).$$
\end{theorem}

We show that the Ball detection function in sample is consistent against general alternatives. Our new detection function can handle the problem of imbalanced sample sizes. As shown in the following theorem, the asymptotic power of the test does not go to zero even if $\eta=\frac{M}{N}$ goes to $0$ or $\infty$.

\begin{theorem}\label{againstall}
	The test based on $V(M,T)/T$ is consistent against any general alternative $H_1$. More specifically,
	$$
	\lim\limits_{(M,T)\rightarrow\infty}Var_{H_1}(V(M,T)/T)=0,
	$$
	and
	$$
	\Lambda:=\liminf\limits_{(M,T)\rightarrow\infty}(E_{H_1}V(M,T)-E_{H_0}V(M,T))/T>0.
	$$
\end{theorem}

\section{Detection of change-points}
\subsection{Hierarchical Algorithm}

Next, we use the Ball detection function in sample to detect change-points in a sequence.
For simplicity, suppose that the sequence $\{Z_i\}_{1\leq i\leq T}$ contains at most one change-point.
The possible change-point location is then estimated by maximizing the detection function:
\begin{equation}\label{singleLoc}
\hat{M}_1=\argmax_{M}V(M,T).
\end{equation}
We use the bootstrap method to estimate the probability that $V(\hat{M}_1,T)$  exceeds a threshold.
If the estimated probability is high enough, $\hat{M}_1$ is the estimated change-point. Otherwise, we proceed as if there does not exist any change-point in the sequence.

It is more complicated if the sequence has multiple change-points. In this case, we estimate the first change-point by
\begin{equation}\label{singleLocLLL}
(\hat{M}_1,\hat{L}_1)=\argmax_{0<M_1<L_1\leq T}V(M_1,L_1).
\end{equation}
From (\ref{singleLocLLL}), we can see that the introduction of $L_1$ here is to alleviate a weakness of bisection algorithm  \citep{matteson2014nonparametric}. Because in each segment, there may exist multiple change-points. If we do not introduce $L_1$,  the value of  $V(\hat{M}_1,T)$ may be lower than $V(\hat{M}_1,\hat{L}_1)$.

Suppose that $k-1$ change-points have been estimated at locations $0<\hat{T}_1<\cdots<\hat{T}_{k-1}<T$, and $\hat{T}_0=0$, $\hat{T}_k=T$. Those change-points partition the sequence into $k$ segments $\mathbf{Z}(\hat{T}_1/\hat{T}_0),\ldots,\mathbf{Z}(\hat{T}_k/\hat{T}_{k-1})$.
In segment $i$, let $C_{ij}^1=\frac{1}{M_i-\hat{T}_{i-1}}\sum_{u=\hat{T}_{i-1}}^{M_i}c(Z_i,Z_j,Z_u),$ $C_{ij}^2=\frac{1}{L_i-M_{i}}\sum_{v=M_{i}}^{L_i}c(Z_i,Z_j,Z_u).$
The Ball detection function in sample of segment $i$ is denoted as
$$
V_i(M_i,L_i)=\frac{(M_i-\hat{T}_{i-1})(L_i-M_i)}{(L_i-\hat{T}_{i-1})^3}\sum_{i=\hat{T}_{I-1}}^{L_i}(C_{ij}^1-C_{ij}^1)^2.
$$
Now let
\begin{equation}\label{kappaformula}
(\hat{M}_{\hat{i}},\hat{L}_{\hat{i}})=\argmax\limits_{1\leq i\leq k-1,\hat{T}_{i-1}<M_i<L_i \leq \hat{T}_i}V_i(M_i,L_i).
\end{equation}
Then  $\hat{M}_{\hat{i}}$ is the $k$-th possible change-point located within segment $\mathbf{Z}(\hat{T}_{\hat{i}}/\hat{T}_{\hat{i}-1}).$
This hierarchical algorithm for estimating multiple change-points is outlined below.
\begin{algorithm}\label{algorithm}
	\caption{Multiple change-points Algorithm}
	\begin{algorithmic}
		\STATE Let the minimum segment size $min=m$, the change-points set $\mathbf{T}=\{0,T\}$.
		\STATE Suppose that $k-1$ change-points have been estimated. This decomposes the observations into $k$ segments.
		\FOR {each $i\in \{1,\ldots,k-1\}$}
		%\STATE
		\IF {the $i-$th segment $\mathbf{Z}(\hat{T}_i/\hat{T}_{i-1})$ is a new segment,}
		\STATE $best_i=0;$
		\FOR{$M_i=\hat{T}_{i-1}+m,\hat{T}_{i-1}+m+1,\ldots,\hat{T}_{i}-m$}
		\FOR{$L_i=\hat{T}_{i-1}+m,\ldots,\hat{T}_{i}$}
		\STATE Compute $V_i(M_i,L_i)$;
		
		\IF{$V_i(M_i,L_i)\geq best$}
		\STATE  $\hat{M}_i=M_i, \hat{L_i}=L_i, V_i(\hat{M}_i,\hat{L}_i)=V_i(M_i,L_i)$;
		\ENDIF
		\ENDFOR
		\ENDFOR
		\ELSE
		\STATE $V_i(\hat{M}_i,\hat{L}_i),\hat{M}_i,\hat{L}_i$ had been calculated.
		\ENDIF
		\ENDFOR
		\STATE $V_i(\hat{M}_{\hat{i}},\hat{L}_{\hat{i}})=\argmax_{0\leq i\leq k-1}V_i(\hat{M}_i,\hat{L}_i).$
		\IF {$V_i(\hat{M}_{\hat{i}},\hat{L}_{\hat{i}})$ exceeds a threshold,}
		\STATE put $\hat{M}_{\hat{i}}$ into $\mathbf{T}$;
		\ELSE
		\STATE there does not exist new change-point.
		\ENDIF
		%\UNTIL{there is no new change-points}
	\end{algorithmic}
\end{algorithm}

\subsection{Hierarchical Significance Testing}

Here, we elaborate the use of the bootstrap method mentioned above.

Theorem \ref{h0} shows that the asymptotic null distribution of  $V(M,T)$ is a mixture of $\chi^2$ distributions. In practice, it is difficult to directly take advantage of the asymptotic null distribution. So, we use the moving block bootstrap \citep{kunsch1989jackknife} to obtain the empirical probabilities.

Given a set of observations $\{Z_t\}_{1\leq t\leq T}$ and the block size $b_T$, we draw a bootstrap resample $\{Z_t^*\}_{1\leq t\leq T}$ as follows: (i) define the $b_T$ dimensional vector $X_t=(Z_t,Z_{t-1},\ldots,Z_{t-b_T+1})$; (ii) resample from block data $\{X_t\}_{1\leq t \leq T-b_T+1}$  with replacement to get pseudo data $\{X_t\}_{1\leq t\leq L}$ which satisfies $T=[Lb_T]$, where $[A]$ denotes the integer part of $A$. Denote the first $T$ elements of $\{X_t\}_{1\leq t \leq L}$ as the bootstrap resample $\{Z_t^*\}_{1\leq t\leq T}$; (iii) repeat steps (i) and (ii) $R$ times. For the $r$-th repetition, denote the maximum value in equation (\ref{singleLocLLL}) based on $\{Z_t^*\}_{1\leq t\leq T}$ by $V(\hat{M}_1^{(r)},\hat{L}_1^{(r)})$; (iv) the approximate probability is estimated by $\{V(\hat{M}_1^{(r)},\hat{L}_1^{(r)}):r=1,\ldots,R\}.$ Denote the threshold of the estimated probability  by $p_T$, if
	 $\frac{\sharp\{V(\hat{M}_1^{(r)},\hat{L}_1^{(r)})\geq V(\hat{M}_1,\hat{L}_1)\}}{R+1}<p_T$,  then $\hat{M}_1$  is a change-point.

In applications, the choice of the block size $b_T$ involves a trade-off. If the block size becomes too small, the moving block bootstrap will destroy the time dependency of the data and the accuracy will deteriorate. But if the block size becomes too large, there will be few blocks to be used. In other words, increasing the block size reduces the bias and captures more persistent dependence, while decreasing the block size reduces the variance as more subsamples are available. Thus, a reasonable trade off is to  consider the mean squared error as the objective criterion to balance the bias and variance. For the linear time series, as proved in  \citet{carlstein1986use}, the value of the block size that minimizes MSE is
$$
b_T^*=\left(\frac{2|\rho|}{1-\rho^2}\right)^{2/3}T^{1/3},
$$
where $\rho$ is the first order autocorrelation.
%\citet{nordman2014convergence} investigate the accuracy of two general non-parametric methods for estimating optimal block lengths for block bootstraps with time series.
Because the construction of MSE depends on the knowledge of the underlying data generating sequence, no optimal result is available in  general.
In this paper, we follow \citet{hong2017testing} and \citet{xiao2007testing} to choose $b_T=\max\{q_T,\bar{q}_T\}$, where
\begin{equation}\label{bp}
q_T=\min\left\{\left[\left(\frac{3T}{2}\right)^{1/3}\left(\frac{2\hat{\rho}}{1-\hat{\rho}^2}\right)^{2/3}\right],\left[8\left(\frac{T}{100}\right)^{1/3}\right]\right\},
\end{equation}
where $\hat{\rho}$ is the estimator of the first autocorrelation of $Z_t$,
$\bar{q}_T$ is the same as (\ref{bp}) except replacing $\hat{\rho}$ with the estimated first order autocorrelation of $Z_t^2$. So the choice of $b_T$ considers the linear dependence and non-linear dependence.

\subsection{Consistency}
The next theorem shows the consistency of the estimated change-point locations under the following assumption.
\begin{assumption}\label{ass1}
Suppose that $\mathbf{Z}(T/0)$  is an absolutely regular sequence which is comprised of two ergodic stationary sequences. Let $\alpha_1 \in (0,1)$ denote the fraction of the observations, such that $\mathbf{Z}(\lfloor\alpha_1 T\rfloor/0)$ be an ergodic stationary sequence with marginal probability measure $\mu_1$,  $\mathbf{Z}(T/\lfloor\alpha_1 T\rfloor)$ the second ergodic stationary sequence with marginal distribution $\mu_2$.  Finally, let $\delta_T$ be a sequence of positive numbers, such that $\delta_T \rightarrow 0$ and $T\delta_T\rightarrow \infty$ as $T\rightarrow \infty$.
\end{assumption}

\begin{theorem}\label{asymptotic1}
	Suppose Assumption \ref{ass1} holds. Let $\hat{M}_1$ be the estimated change-point location from Equation (\ref{singleLoc}) for a sample of size $T$. For all $\epsilon>0$ and $T$ large enough such that $\alpha_1\in[\delta_T,1-\delta_T]$, we have
	\begin{equation*}
	P(\lim_{T\rightarrow\infty}|\frac{\hat{M}_1}{T}-\alpha_1|<\epsilon)=1.
	\end{equation*}
\end{theorem}
This theorem shows that the consistency only requires the size of each segment  increases to $\infty$, but not necessarily at the same rate. Under the Assumption \ref{ass1}, $\alpha_1$ can be close to 0 or 1 when $T\rightarrow\infty$, which is an imbalanced case.

In the multiple change-points situation, we have the following Assumption.
\begin{assumption}\label{ass2mul}
	Suppose that $\{Z_i\}_{1 \leq i \leq T}$ is an absolutely regular sequence. Let $0=T_0<T_1<\ldots<T_k<T_{k+1}=T$, and $\min\limits_{i=1,\ldots,k}|T_i-T_{i-1}|\geq aT^{b}$, with $a>0$ and $0<b\leq1$. For $i=0,1,\ldots,k$,  $\mathbf{Z}(T_i/T_{i-1})$ is an ergodic stationary sequence with marginal probability measure $\mu_i$ and $\mu_i\neq \mu_{i+1}$. \label{key}
	Furthermore, let $\delta_T$ be a sequence of positive numbers, such that $\delta_T \rightarrow 0$ and $T\delta_T\rightarrow \infty$ as $T\rightarrow \infty$.
\end{assumption}

 It is worth noting that we do not assume the upper bounds on the number of change-points $k$, but by specifying  the minimum sample size in each segment. In other words, under Assumption \ref{ass2mul}, as $T\rightarrow\infty$, we can have $k\rightarrow\infty$ change-points.

Analysis of multiple change points can be reduced to the analysis of only two change points under Assumption \ref{ass2mul}.  Let $\alpha_i=T_i/T$, for any  $i\in\{1,\ldots,k\}$. The observations $\mathbf{Z}(T\alpha_i/0)$ can be seen as a random sample from a mixture of probability measures $\{\mu_j:j\leq i\}$, denoted as $\boldsymbol{\mu}_i$. Similarly, observations $\mathbf{Z}(T/T\alpha_{i+1})$ are a sample from a mixture of probability measures $\{\mu_j:j\geq i+1\}$, denoted here as $\boldsymbol{\nu}_i$. The remaining observations are distributed according to some probability measure $\boldsymbol{\xi}_i$. Furthermore, $\boldsymbol{\mu}_i\neq \boldsymbol{\xi}_i$ and $\boldsymbol{\nu}_i\neq \boldsymbol{\xi}_i$.  If one of the previous two inequalities does not hold, we refer to the single change point setting.

Consider any $\alpha$ such that, $\alpha_i\leq\alpha\leq\alpha_{i+1}$. Then, this choice of $\alpha$ will create two mixture probability measures. One with component probability measures $\boldsymbol{\mu}_i$ and $\boldsymbol{\xi}_i$, and the other with component probability measures $\boldsymbol{\nu}_i$ and $\boldsymbol{\xi}_i$. Then, the Ball detection function between these two mixture probability measures is equal to
\begin{equation}\label{multiCP_location}
\begin{split}
&V(\alpha;\frac{\alpha_i}{\alpha}\boldsymbol{\mu}_i+\frac{\alpha-\alpha_i}{\alpha}\boldsymbol{\xi}_i,\frac{1-\alpha_{i+1}}{1-\alpha}\boldsymbol{\nu}_i+\frac{\alpha_{i+1}-\alpha}{1-\alpha}\boldsymbol{\xi}_i)\\=&\alpha(1-\alpha)\iint_{\mathrm{A}\times \mathrm{A}} [\frac{\alpha_i}{\alpha}\boldsymbol{\mu}_i+\frac{\alpha-\alpha_i}{\alpha}\boldsymbol{\xi}_i-\frac{1-\alpha_{i+1}}{1-\alpha}\boldsymbol{\nu}_i-\frac{\alpha_{i+1}-\alpha}{1-\alpha}\boldsymbol{\xi}_i]^{2}(\bar{B}(u,\rho(u,v))\omega_{\alpha} (du)\omega_{\alpha} (dv).
\end{split}
\end{equation}

\begin{theorem}\label{lemma_multiCP_max}
Suppose that Assumption \ref{ass2mul} holds, then the Ball detection function  in equation (\ref{multiCP_location}) is maximized when either $\alpha=\alpha_i$ or $\alpha=\alpha_{i+1}$.
\end{theorem}

By Theorem \ref{lemma_multiCP_max}, $f_i(\alpha)$ is maximized when $\alpha=\alpha_i$ or $\alpha=\alpha_{i+1}$ for $i=1,\ldots,k-1$. Additionally,
define
$$V(\alpha)=\sum_{i=0}^kV(\alpha;\frac{\alpha_i}{\alpha}\boldsymbol{\mu}_i+\frac{\alpha-\alpha_i}{\alpha}\boldsymbol{\xi}_i,\frac{1-\alpha_{i+1}}{1-\alpha}\boldsymbol{\nu}_i+\frac{\alpha_{i+1}-\alpha}{1-\alpha}\boldsymbol{\xi}_i)I(\alpha_i\leq\alpha\leq\alpha_{i+1}).$$
Let $\mathscr{A}_T=\{y\in[\delta_T,1-\delta_T]:V(y)\geq V(\alpha),\forall\alpha\}$. Let $d(x,\mathscr{A}_T)=inf\{|x-y|:y\in\mathscr{A}_T\}$. Then, we have the following Theorem.
\begin{theorem}\label{location_all}
Suppose that Assumption \ref{ass2mul}, and $x\in\mathbb{R}$,  Let $\hat{M}_1$ be the estimated change point as defined by equation (\ref{singleLocLLL}). Then $d(\hat{M}_1/T,\mathscr{A}_T)\xrightarrow{a.s.}0$ as $T\rightarrow\infty$.
\end{theorem}

Repeated applications of Theorem \ref{location_all} can show that as $T\rightarrow\infty$, the first $k$ estimated change points will converge to the true change point locations in the manner described above. With a fixed  threshold of the estimated probability $p_T$, all of the change-points will be estimated.
However, with probability approaching 1 as the sample size increases, the number of change-points determined in this way will be more than the true number of change-points, since any given nominal level of significance implies a nonzero probability of rejecting the null hypothesis when it holds. The hierarchical procedure could be made consistent by adopting a threshold for the test that decrease to zero, at a suitable rate, as the sample size increases\citep{bai1998estimating}.
This is illustrated by the following theorem.
\begin{theorem}
Let $\hat{k}$ be the number of change-points obtained using the hierarchical method based on the statistic (\ref{kappaformula}) applied with threshold $p_T$, and $k$ be the true number of change-points. If $\lim\limits_{T\rightarrow\infty}p_T\rightarrow0$, then under Assumption  \ref{ass2mul}, $\lim\limits_{T\rightarrow\infty}P(\hat{k}=k)=1$.
\end{theorem}
Although the hierarchical algorithm tends to estimate more change-points asymptotically when $p_T$ is fixed, this has little effect in practice. For example, the asymptotic probability of selecting $(k+j)$ change-points, is given by $p_T^j(1-p_T),$ which decreases rapidly.
Furthermore, if there is no change point, that is $k=0$, the probability of selecting at least one change point in our algorithm is
$$
\sum_{j=1}^{\infty}p_T^j(1-p_T)=p_T.
$$ Hence the total rate of type I errors is still $p_T$. This is a distinct feature of our hierarchical procedure because controlling for type I errors is a challenging issue in multiple testings.

\section{Simulation studies}

In this section, we present the numerical performance of the proposed method (BDCP) with $p_T=0.05$ and compare it with several typical methods, including Bayesian method (BCP) \citep{barry1993bayesian}, WBS method  \citep{fryzlewicz2014wild}, the graph-based method-gSeg  \citep{chen2015graph,chu2019asymptotic} and energy distance based method (ECP)  \citep{matteson2014nonparametric}. BCP, WBS, ECP and BDCP can estimate the number of change-points automatically while gSeg can detect only one change-point or an interval.

There are four commonly used criteria for the performance of those methods: the adjusted Rand index, the over segmentation error, the under segmentation error and the Hausdorff distance.
Suppose that the true change-points set is $\mathbf{T}=\{0,T_1,\ldots,T_k,T\}$ and estimated change-points set is $\hat{\mathbf{T}}=\{0,\hat{T}_1,\ldots,\hat{T}_{\hat{k}},T\}$.
Then denote the true segments of series
$\{Z_t\}_{1\leq t\leq T}$ by $\mathbf{Z}=\{\mathbf{Z}(T_1/0),\ldots,\mathbf{Z}(T/T_k)\}$ and the estimated segments by\\
$\mathbf{\hat{Z}}=\{\mathbf{Z}(\hat{T}_1/0),\ldots,\mathbf{Z}(T/\hat{T}_{\hat{k}})\}$.
Consider the pairs of observations that fall into one of the following two sets:\\ $\{S_1\}=$ \{pairs of observations in the same segments under $\mathbf{Z}$ and in same segments under $\mathbf{\hat{Z}}$\};\\ $\{S_2\}=$ \{pairs of observations in different segments under $\mathbf{Z}$ and in different segments under $\mathbf{\hat{Z}}$\}. Denote $\sharp S_1$ and $\sharp S_2$ as the number of pairs of observations in each of these two sets. The Rand index $RI$ is defined as
$$RI=\frac{\sharp S_1+\sharp S_2}{\binom{T}{2}}.$$

%One disadvantage of the Rand index is that it does not {\color{red} take on a constant value (reword here)} when comparing two random clustering \citep{hubert1985comparing}.
Adjusted Rand index $ARI$ is the corrected-for-chance version of the Rand index which is defined as
$$ARI = \frac{RI - E(RI)}{1-E(RI)},$$
in which 1 corresponds to the maximum Rand index value.

On the other hand, we also calculate the distance between $\mathbf{T}$ and $\hat{\mathbf{T}}$ by
$$\zeta(\hat{\mathbf{T}}||\mathbf{T})=\sup_{b\in\mathbf{T}}\inf_{a\in\hat{\mathbf{T}}}|a-b|\ \ and \ \ \zeta(\mathbf{T}||\hat{\mathbf{T}})=\sup_{b\in\hat{\mathbf{T}}}\inf_{a\in\mathbf{T}}|a-b|,$$
which quantify the over-segmentation error and the under-segmentation error, respectively \citep{boysen2009consistencies,zou2014nonparametric}.  The Hausdorff distance  \citep{harchaoui2010multiple} between $\mathbf{T}$ and $\hat{\mathbf{T}}$ is defined as
$$\Delta(\mathbf{T},\hat{\mathbf{T}})=sup\{\zeta(\hat{\mathbf{T}}||\mathbf{T}),\zeta(\mathbf{T}||\hat{\mathbf{T}})\}.$$
Here, we only report the results based on adjusted Rand index. The results under other criteria are deferred to the supplementary material.

Three scenarios are used for comparisons: univariate sequence, multivariate sequence and manifold sequence.
In each scenario, we consider two types of examples, one without change-point, and one with two change-points as follow:
$$\{X_1,X_2,\ldots,X_n,Y_1,Y_2,\ldots,Y_m,X_{n+1},X_{n+2},\ldots,X_{2n}\}.$$
The sample sizes are set to be $n=40$, $m=40,60,80$.
We will repeat each model 400 times and the threshold is at 0.05. To save space, some results of univariate sequences are available on the supplementary material.

\subsection{Multivariate sequence}
In this subsection, we consider the $d=3$ dimensional sequences. Examples 4.1.1-4.1.7 are the sequences with no change-point and Examples 4.1.8-4.1.15 are the models with two change-points.
\begin{itemize}
	\item \textbf{Examples 4.1.1-4.1.3:}
	$$X_t=\epsilon_{t},$$
	$\epsilon_{t}\sim N(0,I_3)$ for Example 4.1.1, $\epsilon_{t}\sim t_3(0,I_3)$ for Example 4.1.2 and $\epsilon_{t}\sim Cauchy(0,I_3)$ for Example 4.1.3.
	
	\item \textbf{Examples 4.1.4-4.1.5}:
	$$X_t=0.5\epsilon_{t}+0.5\epsilon_{t-1},$$
	$\epsilon_{t}\sim N(0,I_3)$ for Example 4.1.4 and $\epsilon_{t}\sim t_3(0,I_3)$ for Example 4.1.5.
	
	\item \textbf{Examples 4.1.6-4.1.7}:
	$$X_{t}=\sigma_{X,t|t-1}\epsilon_{t},$$
	$$\sigma_{X,t|t-1}^{2}=0.02+0.02\sigma_{X,t-1|t-2}^{2}+0.05X_{t}^{2},$$
	$\epsilon_{t}\sim N(0,I_3)$ for Example 4.1.6 and $\epsilon_{t}\sim t_3(0,I_3)$ for Example 4.1.7.
\end{itemize}

\begin{itemize}
	\item \textbf{Example 4.1.8}:
	$$X_t=0.5\epsilon_{t}+0.5\epsilon_{t-1},$$
	$$Y_t=\mu+0.5\epsilon_{t}+0.5\epsilon_{t-1},$$
	$$\mu=(4,4,4),(6,6,6),(8,8,8),,\epsilon_{t}\sim N(0,I_3).$$
	
	\item \textbf{Example 4.1.9}:
	$$X_t=0.5\epsilon_{t}+0.5\epsilon_{t-1},$$
	$$Y_t=\mu+0.5\epsilon_{t}+0.5\epsilon_{t-1},$$
	$$\mu=(4,4,4),(6,6,6),(8,8,8),\epsilon_{t}\sim t_3(0,I_3).$$
	
	\item \textbf{Example 4.1.10}:
	$$X_t=\epsilon_{t},$$
	$$Y_t=\mu+\epsilon_{t},$$
	$$\mu=(4,4,4),(6,6,6),(8,8,8),\epsilon_{t}\sim Cauchy(0,I_3).$$
	
	\item \textbf{Example 4.1.11}:
	$$X_t=0.5\epsilon_{t}+0.5\epsilon_{t-1},$$
	$$Y_t=0.5\epsilon_{t}\sigma+0.5\epsilon_{t-1}\sigma,$$
	$$\sigma=(3,3,3),(5,5,5),(7,7,7),\epsilon_{t}\sim N(0,I_3).$$
	
	\item \textbf{Example 4.1.12}:
	$$X_t=0.5\epsilon_{t}+0.5\epsilon_{t-1},$$
	$$Y_t=0.5\epsilon_{t}\sigma+0.5\epsilon_{t-1}\sigma,$$
	$$\sigma=(3,3,3),(5,5,5),(7,7,7),\epsilon_{t}\sim t_3(0,I_3).$$
	
	\item \textbf{Example 4.1.13}:
	$$X_t=\epsilon_{t},$$
	$$Y_t=\epsilon_{t}\sigma,$$
	$$\sigma=(9,9,9),(16,16,16),(25,25,25),\epsilon_{t}
	\sim Cauchy(0,I_3).$$
	
	\item \textbf{Examples 4.1.14-4.1.15}: \\
	$X_t,Y_t \sim CCC-GARCH(1,1)$\citep{Bollerslev1990Modelling}, let
	$$X_t=\sigma^{X}_{t}\epsilon_t,Y_{t}=\sigma^{Y}_{t}\epsilon_{t},$$
	$$\sigma^{X}_{t}=(\omega_{X}+A^{X}\epsilon_{t-1}+B_{X}\sigma{X}_{t-1}^{2})^{1/2},$$
	$$\sigma^{X}_{t}=(\omega_{Y}+A^{Y}\epsilon_{t-1}+B_{Y}\sigma{Y}_{t-1}^{2})^{1/2}.$$
Let     $\omega^{X}=(0.01,0.01,0.01),$ $A^{X}=diag(0.02,0.03,0.01),$ $B^{X}=diag(0.02,0.02,0.05)$, and

	Case 1:$$\omega^{Y}=2\omega^{X},A_{Y}=4A^{X},B^{Y}=5B^{X},$$
	
	Case 2:$$\omega^{Y}=3\omega^{X},A_{Y}=5A^{X} ,B^{Y}=6B^{X}.$$
	
	Case 3:$$\omega^{Y}=4\omega^{X},A_{Y}=6A^{X} ,B_{Y}=7B^{X}.$$
	$\epsilon_{t}\sim N(0,1)$ for Example 4.1.14 and $\epsilon_{t}\sim t(df=3)$ for Example 4.1.15.
\end{itemize}

Table \ref{Sim_adj_mul_type1} reveals that ECP and BDCP can handle the multivariate stationary series well. BCP works well when the distribution is normal but has a lower adjusted Rand index in $t$ distribution and Cauchy distribution. We do not consider the WBS method because it can not handle the multivariate cases. Tables \ref{Sim_adj_mul_mean}-\ref{Sim_adj_mul_gar} present the results of those examples with two change-points. All four methods have excellent performance in the multivariate normal distribution and multivariate $t$ distribution with location shift.
Examples 4.1.11-4.1.13 consider the scale shift case and Examples 4.1.14-4.1.15 are the popular GARCH models which are also the scale shift case. We can see from Tables  \ref{Sim_adj_mul_var}-\ref{Sim_adj_mul_gar} that BDCP has the best performance in almost all the scale shift cases.

\subsection{Manifold-valued sequence}
In this subsection, we report some manifold-valued examples where ECP can not detect the change-points but BDCP works well. Consider the distribution in a unit circle and let
$$P_1 \sim Unif([-\pi/6,\pi/6)\bigcup[11\pi/6,13\pi/6)),$$
$$P_2 \sim Unif([\pi/3,2\pi/3)\bigcup[7\pi/3,8\pi/3)),$$
$$P_3 \sim Unif([5\pi/6,7\pi/6)\bigcup[17\pi/6,19\pi/6)),$$
$$P_4 \sim Unif([4\pi/3,5\pi/3)\bigcup[10\pi/3,11\pi/3)),$$
$$P_5 \sim Unif([0,4\pi)).$$
We calculate the circular distance which is defined by
\begin{equation}\label{circlediatance}
d(X_i,X_j)=min(|X_i-X_j|,2\pi-|X_i-X_j|).
\end{equation}
We simulate four examples with 0,1,2,3 change-points respectively. Let $n=40, m=40,60,80$.
\begin{itemize}
	\item \textbf{Example 4.2.1}:
	$$\{X_1\ldots,X_{3n}\}\sim P_5.$$
	\item \textbf{Example 4.2.2}:
	$$\{X_1\ldots,X_n\}\sim P_1, \{Y_1\ldots,Y_m\}\sim P_3.$$
	
	\item \textbf{Example 4.2.3}:
	$$\{X_1\ldots,X_n\}\sim P_1, \{Y_1\ldots,Y_m\}\sim P_3,\{X'_1\ldots,X'_n\}\sim P_2.$$
	
	\item \textbf{Example 4.2.4}:
	$$\{X_1\ldots,X_n\}\sim P_1, \{Y_1\ldots,Y_m\}\sim P_3,$$
	$$\{X'_1\ldots,X'_n\}\sim P_2, \{Y'_1\ldots,Y'_m\}\sim P_4.$$
\end{itemize}

\begin{comment}
\begin{figure}[htbp]
\centering
\includegraphics[width=1\textwidth]{sim_exp.pdf}
\caption{\label{sim_exp}}
\end{figure}
\end{comment}

Table \ref{Sim_adj_man_type1} reveals that BCP, WBS, ECP all perform well when there is no change-point. However, they do not work when the sequences have change-points (Table \ref{Sim_adj_man_123}). That is because BCP is based on the normal distribution, and WBS is a CUSUM statistics which does not work in a circular distribution. For ECP, that is because the circular distance is not of strong negative type (Theorem 9.1 in \citet{hjorth1998finite}).  The gSeg method can detect change-points when the number of change-points is one or two but do not perform well in Example 4.2.4. BDCP has a remarkable performance in all these examples.

\section{Real data analysis}

\subsection{Wind direction of Yunnan-Guizhou Plateau}

Monsoon is used to describe seasonal changes in atmospheric circulation and precipitation associated with the asymmetric heating of the land and sea. The major monsoon systems in the world consist of West African Monsoon (WAM), Indian summer monsoon (ISM), East Asian Monsoon (EAM) and so on. In this subsection, we analyze the wind direction data of Yunnan-Guizhou Plateau ($105^\circ$E, $27^\circ$N) from 06/01/2015 to 10/30/2015. The data are available in R package \emph{rWind}. Yunnan-Guizhou Plateau is located in southwest China, with local climate influenced by both ISM and EAM  \citep{sirocko1996teleconnections}\citep{li2014synchronous}(Fig. 1A and S1). Strict spatial boundaries between the ISM and the ASM are difficult to define \citep{cheng2012global} though previous researchers have suggested $103^\circ$E as the dividing line on the basis of summer prevailing winds.

\begin{figure}[!htb]
	\centering
	\includegraphics[width=1\textwidth]{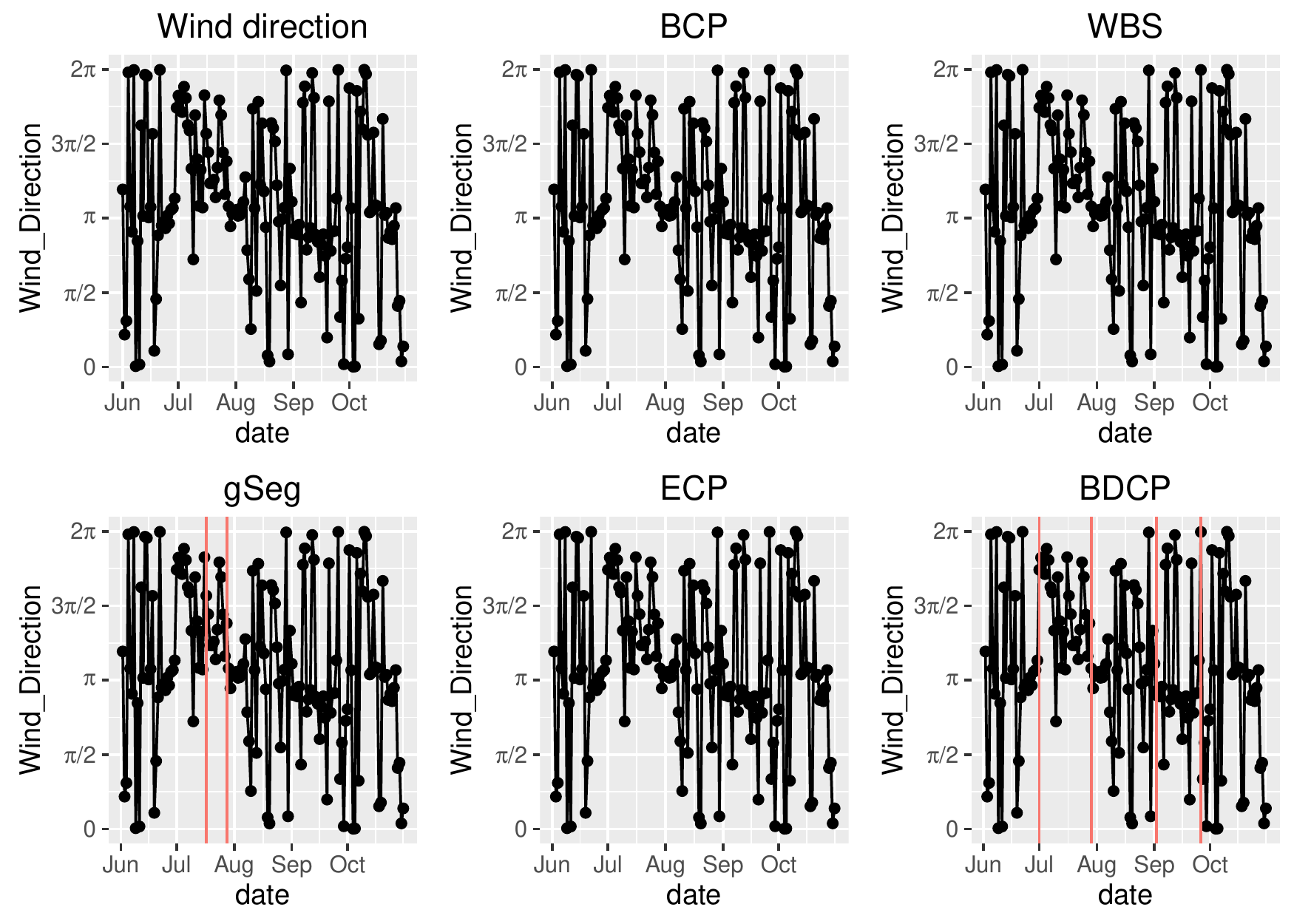}
	\caption{\label{wind1}Wind direction of Yunnan-Guizhou Plateau from 06/01/2015 to 10/30/2015 and the performance of BCP, WBS, gSeg, ECP and BDCP. The y-axis shows wind direction. Degree 0 and $2\pi$ represent due North, $\pi/2$ represents due East, $2\pi$ represents due South and $3\pi/2$ represents due West.}
\end{figure}

Daily wind directions are shown in the top-left of Figure \ref{wind1}. Note that degree 0 represents due North, $\pi/2$ represents due East, $\pi$ represents due South and $3\pi/2$ represents due West.
We can see that the wind directions are distributed in almost all directions. In the beginning, the most widely distributed direction is the southwest wind from the Indian Ocean, and then turns smoothly to southeast, which is from the Pacific Ocean. In particular, Yunnan-Guizhou Plateau was mostly influenced by ISM in June and July. After July, the influence of EAM gradually increased  \citep{li2015introduction}.

To detect the change-point in the wind direction series, we calculate the circular distance between the daily direction as defined in (\ref{circlediatance}).

The performance of the five methods is shown  in Figure \ref{wind1}.
BCP, WBS and ECP can not detect any change-point, as seen in the simulation studies in subsection 4.2. gSeg detects an interval between ``07/17/2015'' and ``07/28/2015''. BDCP estimates four change-points located at ``07/01/2015'', ``07/29/2015'', ``09/02/2015'' and ``09/26/2015''.

%The goal of this work is to investigate the coupled interannual variability of the ISM and ASM, and the dynamic role of intervening high topography in this linkage.

\begin{figure}[!htb]
	\centering
	\includegraphics[width=1\textwidth]{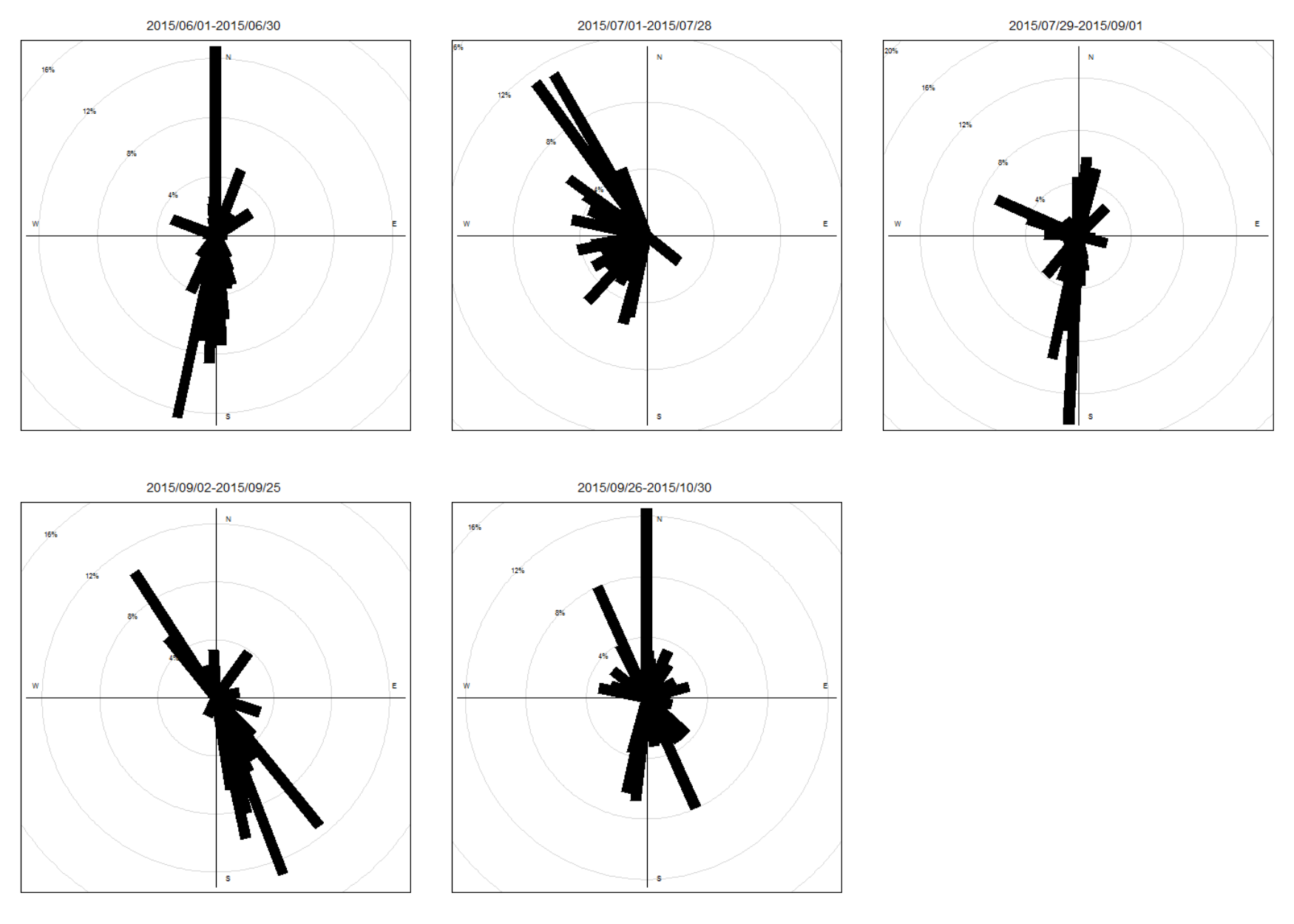}
	\caption{\label{five_periods}The wind rose plot for the five periods detected by BDCP.}
\end{figure}

To visualize the result of BDCP, Figure \ref{five_periods} depicts the wind rose plot for the five periods detected by BDCP. We can see the significant changes of the direction distribution especially between 07/29/2015 - 09/01/2015 and 09/02/2015 - 09/25/2015. The wind directions are almost southwest or west in June and July, then  turn to southeast in September \citep{li2015introduction}. As mentioned in \citet{hillman20178}, 75\% of the average annual precipitation falls in the months of June-September associated with the ISM, and the ISM gets weaker during June and July because isolation decreases by 2-3\%. BDCP can perfectly detect the change of influence between ISM and EAM.

\subsection{Bitcoin price}
Bitcoin is the most popular form of cryptocurrency in recent years.  According to research of Cambridge University in 2017  \citep{hileman2017global}, there are 2.9 to 5.8 million unique cryptocurrency wallet users, most of whom use Bitcoin. One of the known features of Bitcoin is its high volatility.
Bitcoin is not a denominated flat currency and there is no central bank overseeing the issuing of Bitcoin,  its price is thus driven solely by the investors. Using the weekly data over 2010-2013 period,   \citet{briere2015virtual} showed that Bitcoin investment had some high distinctive features, including exceptionally high average return and volatility. Hence, accurately fitting its variation is important \citep{chu2015statistical}.

%Bitcoin is relatively new

\begin{figure}[!htb]
	\centering
	\includegraphics[width=1\textwidth]{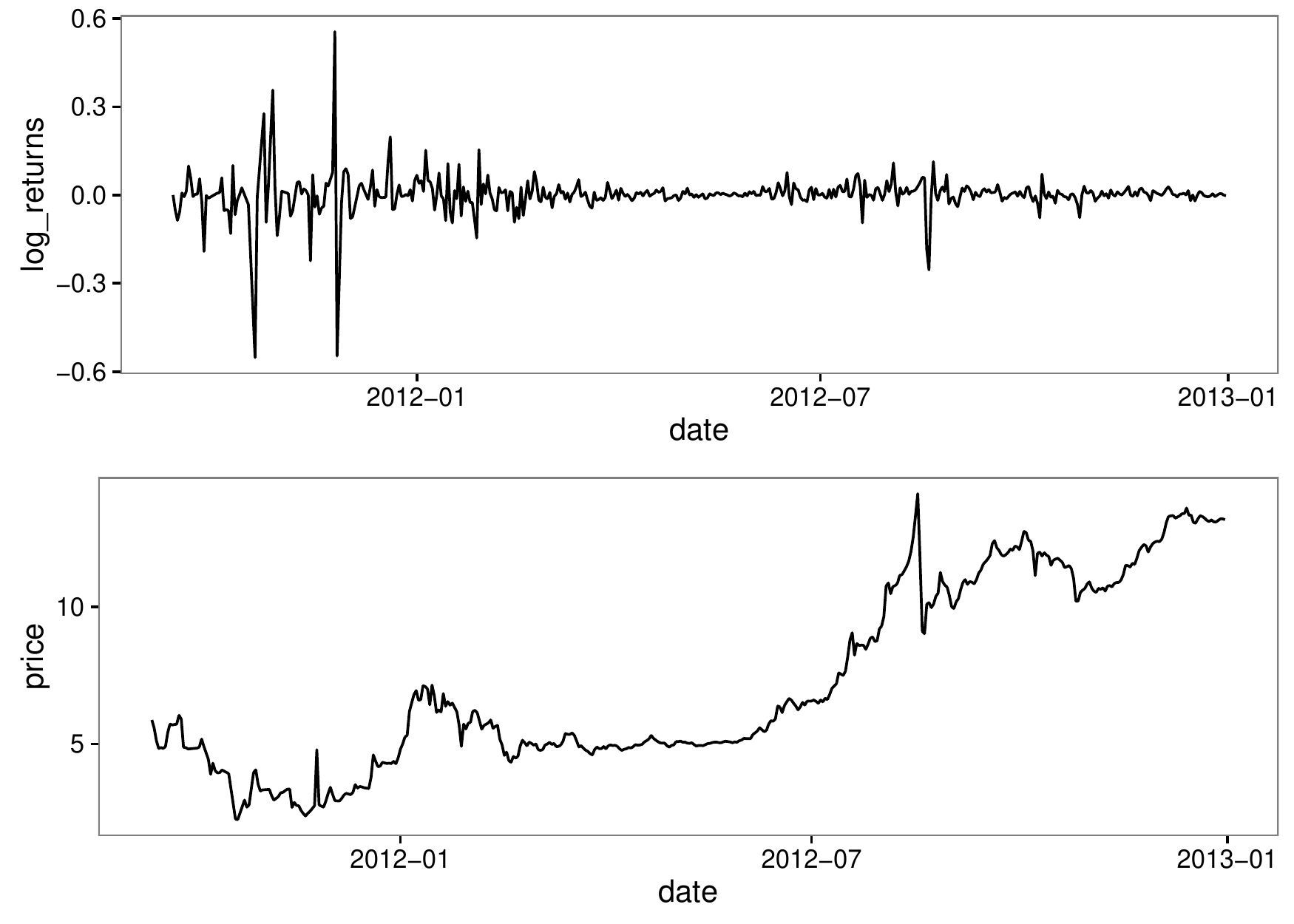}
	\caption{\label{bitcoin} US dollar-Bitcoin exchange rate and daily log-returns from 09/13/2011-12/31/2012.}
\end{figure}

Bitcoin can be exchanged for other currencies, products, and services in legal or black markets. \citet{chu2015statistical} measured the volatility of Bitcoin exchange rate against six major currencies. They found that the behavior of Bitcoin was sharply different from those currencies; its interquartile range was much wider, its skewness was much more negative, its kurtosis was much more peaked and its variance was much larger. Bitcoin showed the highest annualized volatility of percentage change in daily exchange rates. In this subsection, we detect the change-points of daily log-return of Bitcoin using methods, BCP, WBS, gSeg, ECP and BDCP. The datasets are available on \url{http://api.bitcoincharts.com/v1/csv/bitstampUSD.csv.gz}. Figure \ref{bitcoin} displays the  exchange rate of Bitcoin and daily log-returns during 09/13/2011 - 12/31/2012.

\begin{figure}[!htb]
	\centering
	\includegraphics[width=1\textwidth,height=0.9\textheight]{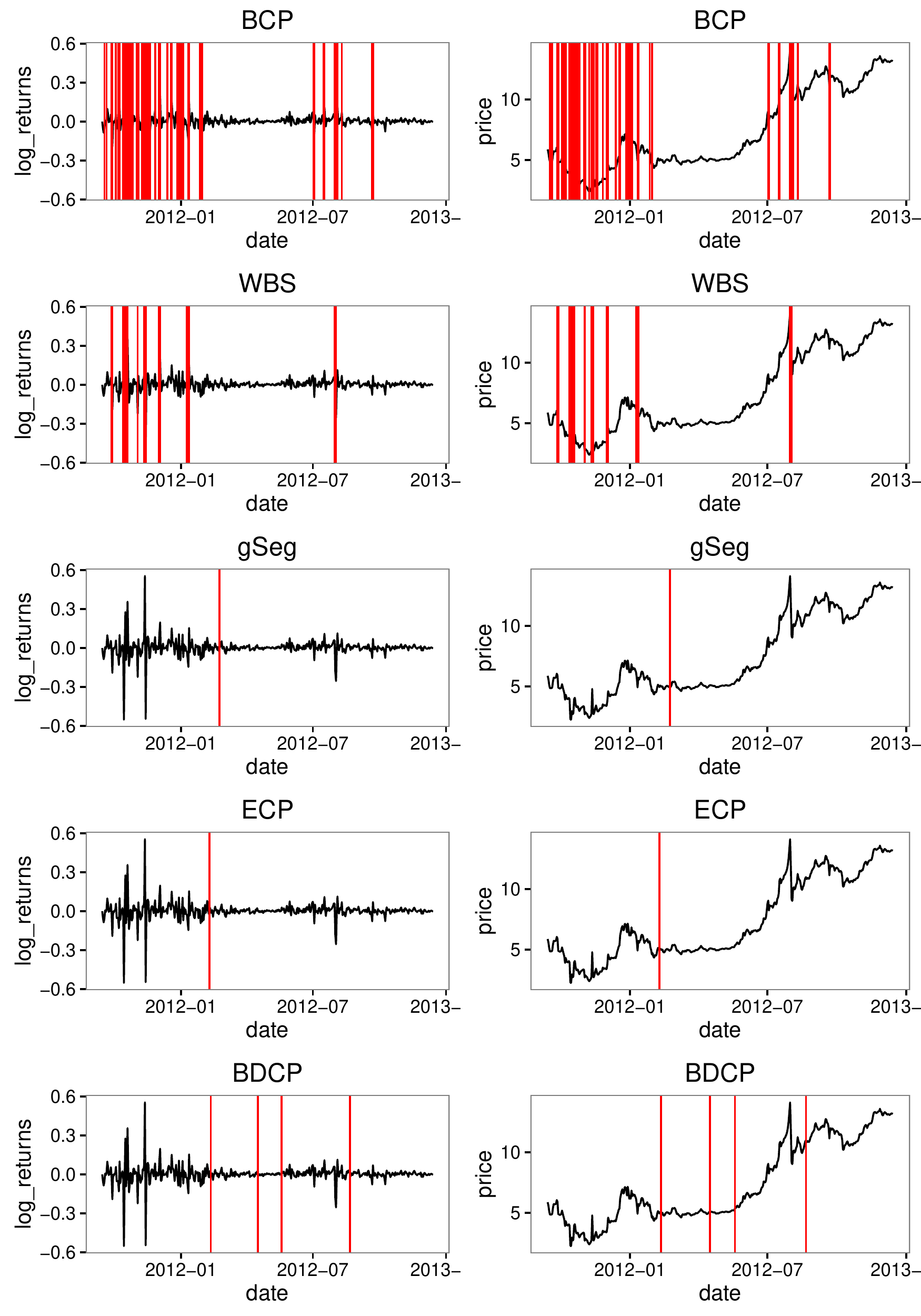}
	\caption{\label{bitcoincp} The results of BCP, WBS, geg, ECP, BDCP on the daily log-returns of Bitcoin series from 09/13/2011-12/31/2012.  The y-axis shows daily log-returns of Bitcoin series. The red lines are the change-point locations detected by BCP, WBS, geg, ECP, BDCP. } 
\end{figure}

Figure \ref{bitcoincp} compares the performance of the five methods. BCP and WBS can not handle the severe volatility at the beginning of the sequence.
gSeg detects one change-point at ``02/23/2012'', and ECP estimates a change-point at ``02/09/2012''.
BDCP detects four change-points at ``02/11/2012'', ``04/16/2012'', ``05/19/2012'', and ``08/21/2012''.

On February 11, 2012, Paxum, an online payment service and popular means for exchanging Bitcoin announced it would cease all dealings related to the currency due to the concerns of its legality. Two days later, regulatory issues surrounding money transmission compelled the popular Bitcoin exchange and service firm TradeHill to terminate its business and immediately began selling its Bitcoin assets to refund its customers and creditors. Bitcoin trading started to cool down during that period.

After May 19, the price of Bitcoin had increased from \$5.07 to the maximum \$14.14 on August 17 and kept at that level after that. The reasons for the rise were many. Lots of online articles on this subject expressed the same message: Bitcoin was now going mainstream. WordPress, ranked by Alexa as the 21st most popular site in the world, started to accept Bitcoin for payment on November, 2012.

The variances of these five stages detected by BDCP are: 0.1054, 0.0182, 0.0059, 0.0458, and 0.0177. The daily log-return sequence was very flat and the price almost did not change during period 04/16/2012 - 05/18/2012. But it was volatile during other periods from Figure \ref{bitcoincp}.

\section{Conclusion}
We developed a change-point detection procedure for weakly dependent sequences. Our key idea lies in the novel measure of Ball detection function.
We proved the asymptotic properties of its sample statistic for absolutely regular sequences.  Extensive simulation studies demonstrated that our method had a superior performance to other existing methods in various settings. Two real data analyses indicated that our method was useful in analyzing non-Euclidean sequences with various change points and led to insightful understanding of the data. Also, our method is robust since our test statistic is rank-based.

We will further investigate Ball detection function and its related concepts. For example,  the current computational complexity of our proposed algorithm is $O(kT^2\log T)$, where $k$ is the number of change-points, and $T$ is the length of the sequence. It will be useful to find an algorithm with a lower computational complexity.

\bibliographystyle{ECA_jasa}
\bibliography{reference}

\newpage
\begin{table}[H] %\scriptsize
	\renewcommand{\arraystretch}{0.7}
	\caption{The performance of adjusted Rand index for multivariate series with no change-point. The highest average adjusted Rand index is highlighted in bold. The last four columns refer to the adjusted Rand index ratio between the four methods to BDCP.}
	{
		\begin{center}
			\begin{tabular}{lr|rrr}
				\hline
				Example & \multicolumn{1}{l|}{BDCP} & \multicolumn{1}{l}{BCP/BDCP} & \multicolumn{1}{l}{gSeg/BDCP} & \multicolumn{1}{l}{ECP/BDCP} \\
				\hline
				4.1.1 & 0.955 & \textbf{1.042} & 0.000 & 0.963 \\
				4.1.2 & \textbf{0.945} & 0.048 & 0.000 & 0.995 \\
				4.1.3 & \textbf{0.930} & 0.005 & 0.000 & 0.984 \\
				4.1.4 & 0.955 & 0.984 & 0.000 & \textbf{1.047} \\
				4.1.5 & 0.825 & 0.248 & 0.000 & \textbf{1.206} \\
				4.1.6 & \textbf{0.965} & 0.969 & 0.000 & 0.979 \\
				4.1.7 & \textbf{0.920} & 0.087 & 0.000 & 0.995 \\
				\hline
			\end{tabular}\label{Sim_adj_mul_type1}
	\end{center}}
\end{table}%

\begin{table}[H] %\scriptsize
	\renewcommand{\arraystretch}{0.7}
	\caption{The performance of adjusted Rand index for Examples 4.1.8 - 4.1.10 with change in mean. The highest average adjusted Rand index is highlighted in bold. The last four columns refer to the adjusted Rand index ratio between the four methods to BDCP.}
	{
		\begin{center}
			\begin{tabular}{cccr|rrr}
				\hline
				Example & m     & $\mu$   & \multicolumn{1}{l|}{BDCP} & \multicolumn{1}{l}{BCP/BDCP} & \multicolumn{1}{l}{gSeg/BDCP} & \multicolumn{1}{l}{ECP/BDCP} \\
				\midrule
				\multirow{9}[6]{*}{4.1.8} & \multirow{3}[2]{*}{40} & 4     & 0.962 & 0.947 & 0.964 & \textbf{1.018} \\
				&       & 6     & 0.996 & 0.999 & 0.995 & \textbf{1.004} \\
				&       & 8     & 0.995 & \textbf{1.005} & \textbf{1.005} & \textbf{1.005} \\
				\cmidrule{2-7}          & \multirow{3}[2]{*}{60} & 4     & 0.970 & 0.941 & 0.969 & \textbf{1.008} \\
				&       & 6     & 0.994 & 0.999 & 1.000 & \textbf{1.006} \\
				&       & 8     & 0.987 & \textbf{1.013} & 1.012 & \textbf{1.013} \\
				\cmidrule{2-7}          & \multirow{3}[2]{*}{80} & 4     & 0.969 & 0.981 & 0.979 & \textbf{1.009} \\
				&       & 6     & 0.991 & 1.005 & 1.005 & \textbf{1.005} \\
				&       & 8     & 0.987 & \textbf{1.013} & \textbf{1.013} & 1.012 \\
				\midrule
				\multirow{9}[6]{*}{4.1.9} & \multirow{3}[2]{*}{40} & 4     & 0.979 & 0.853 & 0.996 & \textbf{1.015} \\
				&       & 6     & 0.987 & 0.902 & 1.004 & \textbf{1.012} \\
				&       & 8     & 0.990 & 0.938 & 1.006 & \textbf{1.010} \\
				\cmidrule{2-7}          & \multirow{3}[2]{*}{60} & 4     & 0.972 & 0.823 & 1.001 & \textbf{1.020} \\
				&       & 6     & 0.973 & 0.898 & 1.016 & \textbf{1.026} \\
				&       & 8     & 0.988 & 0.914 & 1.009 & \textbf{1.009} \\
				\cmidrule{2-7}          & \multirow{3}[2]{*}{80} & 4     & 0.961 & 0.838 & 1.023 & \textbf{1.027} \\
				&       & 6     & 0.971 & 0.884 & \textbf{1.026} & 1.025 \\
				&       & 8     & 0.966 & 0.900 & \textbf{1.033} & 1.032 \\
				\hline
				\multirow{9}[6]{*}{4.1.10} & \multirow{3}[2]{*}{40} & 4     & \textbf{0.979} & 0.444 & 0.993 & 0.836 \\
				&       & 6     & \textbf{0.987} & 0.496 & 0.999 & 0.949 \\
				&       & 8     & \textbf{0.988} & 0.514 & 0.999 & 0.975 \\
				\cmidrule{2-7}          & \multirow{3}[2]{*}{60} & 4     & 0.981 & 0.437 & \textbf{1.000} & 0.848 \\
				&       & 6     & 0.985 & 0.469 & \textbf{1.002} & 0.937 \\
				&       & 8     & 0.986 & 0.506 & \textbf{1.004} & 0.989 \\
				\cmidrule{2-7}          & \multirow{3}[2]{*}{80} & 4     & 0.975 & 0.389 & \textbf{1.012} & 0.842 \\
				&       & 6     & 0.980 & 0.420 & \textbf{1.012} & 0.950 \\
				&       & 8     & 0.983 & 0.476 & \textbf{1.010} & 0.983 \\
				\hline
			\end{tabular}%
			\label{Sim_adj_mul_mean}
	\end{center}}
\end{table}

\begin{table}[H] %\scriptsize
	\renewcommand{\arraystretch}{0.7}
	\caption{The performance of adjusted Rand index for Examples 4.1.11 - 4.1.13 with change in scale. The highest average adjusted Rand index is highlighted in bold. The last four columns refer to the adjusted Rand index ratio between the four methods to BDCP.}
	{
		\begin{center}
			\begin{tabular}{cccr|rrr}
				\hline
				Example & m     & $\sigma$   & \multicolumn{1}{l|}{BDCP} & \multicolumn{1}{l}{BCP/BDCP} & \multicolumn{1}{l}{gSeg/BDCP} & \multicolumn{1}{l}{ECP/BDCP} \\
				\hline
				\multirow{9}[6]{*}{4.1.11} & \multirow{3}[2]{*}{40} & 3     & \textbf{0.786} & 0.384 & 0.888 & 0.052 \\
				&       & 5     & \textbf{0.931} & 0.622 & 0.911 & 0.632 \\
				&       & 7     & \textbf{0.956} & 0.663 & 0.949 & 0.941 \\
				\cmidrule{2-7}          & \multirow{3}[2]{*}{60} & 3     & \textbf{0.834} & 0.206 & 0.787 & 0.036 \\
				&       & 5     & \textbf{0.942} & 0.408 & 0.908 & 0.646 \\
				&       & 7     & \textbf{0.954} & 0.471 & 0.985 & 0.971 \\
				\cmidrule{2-7}          & \multirow{3}[2]{*}{80} & 3     & \textbf{0.822} & 0.155 & 0.787 & 0.052 \\
				&       & 5     & \textbf{0.941} & 0.248 & 0.919 & 0.624 \\
				&       & 7     & \textbf{0.955} & 0.272 & 0.978 & 0.981 \\
				\hline
				\multirow{9}[6]{*}{4.1.12} & \multirow{3}[2]{*}{40} & 3     & \textbf{0.521} & 0.810 & 0.964 & 0.123 \\
				&       & 5     & \textbf{0.838} & 0.621 & 0.885 & 0.443 \\
				&       & 7     & \textbf{0.906} & 0.603 & 0.905 & 0.715 \\
				\cmidrule{2-7}          & \multirow{3}[2]{*}{60} & 3     & \textbf{0.597} & 0.616 & 0.915 & 0.101 \\
				&       & 5     & \textbf{0.833} & 0.475 & 0.870 & 0.459 \\
				&       & 7     & \textbf{0.893} & 0.477 & 0.920 & 0.776 \\
				\cmidrule{2-7}          & \multirow{3}[2]{*}{80} & 3     & \textbf{0.640} & 0.480 & 0.780 & 0.086 \\
				&       & 5     & \textbf{0.842} & 0.359 & 0.809 & 0.469 \\
				&       & 7     & \textbf{0.898} & 0.356 & 0.893 & 0.751 \\
				\hline
				\multirow{9}[6]{*}{4.1.13} & \multirow{3}[2]{*}{40} & 9     & \textbf{0.686} & 0.618 & 0.914 & 0.058 \\
				&       & 16    & \textbf{0.854} & 0.488 & 0.874 & 0.109 \\
				&       & 25    & \textbf{0.907} & 0.492 & 0.883 & 0.141 \\
				\cmidrule{2-7}          & \multirow{3}[2]{*}{60} & 9     & \textbf{0.707} & 0.497 & 0.932 & 0.054 \\
				&       & 16    & \textbf{0.885} & 0.424 & 0.884 & 0.089 \\
				&       & 25    & \textbf{0.918} & 0.406 & 0.908 & 0.120 \\
				\cmidrule{2-7}          & \multirow{3}[2]{*}{80} & 9     & \textbf{0.736} & 0.443 & 0.841 & 0.067 \\
				&       & 16    & \textbf{0.883} & 0.356 & 0.900 & 0.053 \\
				&       & 25    & \textbf{0.927} & 0.350 & 0.924 & 0.061 \\
				\hline
			\end{tabular}%
			\label{Sim_adj_mul_var}
	\end{center}}
\end{table}

\begin{table}[H]
	\renewcommand{\arraystretch}{0.7}
	\caption{The performance of adjusted Rand index for Examples 4.1.14 and 4.1.15 with change in parameters of GARCH model. The highest average adjusted Rand index is highlighted in bold. The last four columns refer to the adjusted Rand index ratio between the four methods to BDCP.}
	{
		\begin{center}
			\begin{tabular}{cccr|rrr}
				\hline
				Example & m     & case & \multicolumn{1}{l|}{BDCP} & \multicolumn{1}{l}{BCP/BDCP} & \multicolumn{1}{l}{gSeg/BDCP} & \multicolumn{1}{l}{ECP/BDCP} \\
				\hline
				\multirow{9}[6]{*}{4.1.14} & \multirow{3}[2]{*}{40} & 1     & \textbf{0.786} & 0.384 & 0.771 & 0.052 \\
				&       & 2     & \textbf{0.931} & 0.622 & 0.875 & 0.632 \\
				&       & 3     & \textbf{0.956} & 0.663 & 0.925 & 0.941 \\
				\cmidrule{2-7}          & \multirow{3}[2]{*}{60} & 1     & \textbf{0.834} & 0.206 & 0.689 & 0.036 \\
				&       & 2     & \textbf{0.942} & 0.408 & 0.827 & 0.646 \\
				&       & 3     & \textbf{0.954} & 0.471 & 0.948 & 0.971 \\
				\cmidrule{2-7}          & \multirow{3}[2]{*}{80} & 1     & \textbf{0.822} & 0.155 & 0.619 & 0.052 \\
				&       & 2     & \textbf{0.941} & 0.248 & 0.811 & 0.624 \\
				&       & 3     & \textbf{0.955} & 0.272 & 0.938 & 0.981 \\
				\hline
				\multirow{9}[6]{*}{4.1.15} & \multirow{3}[2]{*}{40} & 1     & 0.521 & 0.810 & \textbf{1.086} & 0.123 \\
				&       & 2     & \textbf{0.838} & 0.621 & 0.834 & 0.443 \\
				&       & 3     & \textbf{0.906} & 0.603 & 0.877 & 0.715 \\
				\cmidrule{2-7}          & \multirow{3}[2]{*}{60} & 1     & \textbf{0.597} & 0.616 & 0.859 & 0.101 \\
				&       & 2     & \textbf{0.833} & 0.475 & 0.801 & 0.459 \\
				&       & 3     & \textbf{0.893} & 0.477 & 0.897 & 0.776 \\
				\cmidrule{2-7}          & \multirow{3}[2]{*}{80} & 1     & \textbf{0.640} & 0.480 & 0.705 & 0.086 \\
				&       & 2     & \textbf{0.842} & 0.359 & 0.787 & 0.469 \\
				&       & 3     & \textbf{0.898} & 0.356 & 0.880 & 0.751 \\
				\hline
			\end{tabular}\label{Sim_adj_mul_gar}
	\end{center}}
\end{table}

\begin{table}[H]
	\renewcommand{\arraystretch}{0.7}
	\caption{The performance of adjusted Rand index for manifold series with no change-point. The highest average adjusted Rand index is highlighted in bold. The last four columns refer to the adjusted Rand index ratio between the four methods to BDCP.}
	{
		\begin{center}
			\begin{tabular}{lrr|rrrr}
				\hline
				Example & \multicolumn{1}{l}{T} & \multicolumn{1}{l|}{BDCP} & \multicolumn{1}{l}{BCP/BDCP} & \multicolumn{1}{l}{WBS/BDCP} & \multicolumn{1}{l}{gSeg/BDCP} & \multicolumn{1}{l}{ECP/BDCP} \\
				\hline
				4.2.1 & 120   & 0.915 & \textbf{1.087} & 1.011 & 0.000 & 1.027 \\
				4.2.1 & 140   & 0.945 & \textbf{1.048} & 1.005 & 0.000 & 1.011 \\
				4.2.1 & 160   & 0.940 & \textbf{1.064} & 1.016 & 0.000 & 1.005 \\
				\hline
			\end{tabular}\label{Sim_adj_man_type1}
	\end{center}}
\end{table}%

\begin{table}[H]
	\renewcommand{\arraystretch}{0.7}
	\caption{The performance of adjusted Rand index for manifold series with 1,2,3 change-points. The highest average adjusted Rand index is highlighted in bold. The last four columns refer to the adjusted Rand index ratio between the four methods to BDCP.}
	{
		\begin{center}
			\begin{tabular}{rrr|rrrr}
				\hline
				\multicolumn{1}{l}{Example} & \multicolumn{1}{l}{m} & \multicolumn{1}{l|}{BDCP} & \multicolumn{1}{l}{BCP/BDCP} & \multicolumn{1}{l}{WBS/BDCP} & \multicolumn{1}{l}{gSeg/BDCP} & \multicolumn{1}{l}{ECP/BDCP} \\
				\hline
				\multicolumn{1}{r}{\multirow{3}[2]{*}{4.2.2}} & 40    & \textbf{0.989} & 0.007 & 0.029 & 0.961 & 0.123 \\
				& 60    & \textbf{0.986} & 0.012 & 0.016 & 0.972 & 0.074 \\
				& 80    & \textbf{0.985} & 0.003 & 0.019 & 0.980 & 0.037 \\
				\hline
				\multicolumn{1}{r}{\multirow{3}[2]{*}{4.2.3}} & 40    & \textbf{0.994} & 0.004 & 0.054 & 0.901 & 0.053 \\
				& 60    & 0.993 & 0.007 & 0.088 & \textbf{1.004} & 0.082 \\
				& 80    & 0.983 & 0.008 & 0.123 & \textbf{1.017} & 0.071 \\
				\hline
				\multicolumn{1}{r}{\multirow{3}[2]{*}{4.2.4}} & 40    & \textbf{0.995} & 0.003 & 0.317 & 0.714 & 0.096 \\
				& 60    & \textbf{0.994} & 0.005 & 0.414 & 0.730 & 0.114 \\
				& 80    & \textbf{0.994} & 0.009 & 0.482 & 0.756 & 0.110 \\
				\hline
			\end{tabular}\label{Sim_adj_man_123}
	\end{center}}
\end{table}%

%\section{Appendix}
%\appendix{Supplement A}
%\sname{Supplement A}\label{suppA}
%\stitle{Title of the Supplement A}

%\end{supplement}

\end{document}